\documentclass{conm-p-l}

\usepackage{amsmath,  amssymb,amsbsy,amsfonts,amsthm,latexsym,amsopn,amstext, amsxtra,euscript,amscd}
\usepackage{hyperref}   

\usepackage{graphicx}

\hypersetup{
  colorlinks   = true, 
  urlcolor     = blue, 
  linkcolor    = blue, 
  citecolor   = red 
}

\usepackage{amsrefs}

\newtheorem{thm}{Theorem}
\newtheorem{lem}{Lemma}

\newtheorem{prop}{Proposition}

\newtheorem{exa}{Example}
\newtheorem{rem}{Remark}

\newtheorem{prob}{Problem}

\numberwithin{equation}{section}

\def\R{\mathbb R}
\def\F{\mathbb F}
\def\Z{\mathbb Z}
\def\C{\mathbb C}
\def\P{\mathbb P}
\def\Q{\mathbb Q}

\def\X{\mathcal X}

\def\p{\bar p}
\def\gl{\mbox{GL}}
\def\g{\gamma}

\begin{document}

\title{Self-inversive polynomials,   curves, and codes}


\author{D. Joyner}
\address{Department of Mathematics, US Naval Academy, Annapolis, MD, 21402}
\email{wdjoyner@gmail.com}

\author{T. Shaska}
\address{546 Mathematics and Science Center, Rochester, MI, 48309}
\email{shaska@oakland.edu}


\date{}

\begin{abstract}
We study connections between  self-inversive and self-reciprocal polynomials, reduction theory of binary forms, minimal models of curves, and formally self-dual codes.
We prove that if $\X$ is a superelliptic curve defined over $\C$ and its reduced automorphism group is nontrivial or not isomorphic to a cyclic group, then we can write   its equation as $y^n = f(x)$ or $y^n = x f(x)$, where $f(x)$ is a self-inversive or self-reciprocal polynomial. 
Moreover, we state a conjecture on the coefficients of the   zeta polynomial of extremal formally self-dual codes. 
\end{abstract}

\maketitle



\def\aa{\bar a}
\def\z{\bar z}
\def\H{\mathcal H}
\def\img{\mbox{Img }}

\section{Introduction}

Self-inversive and self-reciprocal polynomials have been studied extensively in the last few decades due to their connections to complex functions and number theory.  In this paper we explore the connections between such polynomials to algebraic curves, reduction theory of binary forms, and coding theory.  While connections to coding theory have been explored by many authors before we are not aware of any previous work that explores the connections of self-inversive and self-reciprocal polynomials to superelliptic curves and reduction theory.  




In section \ref{sect-1}, we give a geometric introduction to inversive
and 
reciprocal polynomials of a given polynomial.  We motivate such
definitions via the 
transformations of the complex plane which is the original motivation
to study such 
polynomials.  It is unclear who coined the names inversive,
reciprocal, palindromic, 
and antipalindromic, but it is obvious that inversive  come from the
inversion 
$z \mapsto \frac 1 {\bar z}$ 
and reciprocal from the reciprocal map  $z \mapsto \frac 1 { z}$ of the complex plane. 

We take the point of view of the reduction theory of binary forms.
While this is an 
elegant and beautiful theory for binary quadratics, it is rather
technical for higher 
degree forms. However,  the inversion plays an important role on
reduction as can 
be seen from section 2 and from \cite{b-sh} and \cite{beshaj}. We are
not aware of 
other authors have explored the connection between reduction theory
and 
self-inversive and self-reciprocal polynomials before even though the overlap is quite obvious. 

We state some of the main results of self-inversive polynomials
including the {\it middle coefficient conjecture} (\ref{eqn:MCC}) and results on the location 
of the roots of such polynomials. 
Self-inversive polynomials over 
$\Q$, $\R$, and $\C$ are discussed and a few recent results on the
height of such polynomials.   The normal references here are 
\cite{conrad, lakatos, ohara, vieira, schinzel-1, schinzel-2, joyner, L1, L2}.
Further,  we discuss the roots of the self-inversive polynomials.
There is a huge amount of literature on this topic including several
conjectures. It is the location of such roots that makes
self-inversive polynomials 
interesting in reduction theory, coding theory, and other areas of mathematics. 
 An attempt at a converse to this
conjecture is discussed in \S \ref{sec:SRPs}.

In section \ref{sect-2} it is given an account of how self-inversive polynomials can be used to determine minimal polynomials of superelliptic curves with extra automorphisms.  This is a new idea spurred by Beshaj's thesis \cite{beshaj} and \cite{reduction}
and has some interesting relations between two different areas of mathematics, namely the theory of algebraic curves and the theory of self-inversive polynomials.  Further details in this direction are planned in \cite{b-sh}. In  this section we prove that for any superelliptic curve with reduced automorphism group not trivial and not isomorphic to a cyclic group we can write the equation of the curve as 
$ y^n = f(x)$   or $y^n = x f(x)$, where $f(x)$ is a palindromic, antipalindromic, or self-inversive polynomial.  Indeed, we can say more since in each case when the automorphism group of the curve we can determine the polynomial $f(x)$ specifically. 

In section \ref{sect-3} we explore connections of self-inversive and self-reciprocal polynomials to reduction theory of binary forms.  We show that
self-inversive polynomials which have all roots on the unit circle correspond to the totally real forms.  The reduction theory for such forms is simpler than for other forms since the Julia quadratic of any degree $n$ form $f(x, y)$ is a factor of a degree $(n-1)(n-2)$ covariant $G_f (x, y)$ given in terms of the partial derivatives of $f$; see \cite{reduction}. We prove that for $f$ palindromic, $G_f$ is self-inversive and if $f$ is palindromic of odd degree then $G_f$ is palindromic. Moreover, we determine explicitly which self-inversive polynomials $f$ with all roots on the unit circle are reduced.

In section \ref{sect-4} we discuss the Riemann hypothesis for formal
weight enumerators of codes and its relation to the self-inversive
polynomials.  We state several open problems which relate to Riemann
hypotheses for extremal formal weight enumerators of codes. 

Most of the results obtained here, with the necessary adjustments,  can be extended to curves defined over fields of positive characteristic.  In \cite{Sa2} equations of superelliptic curves are also determined over such fields. The main question that comes from the connection between self-inversive and self-reciprocal polynomials and reduction theory is whether such polynomials are actually reduced.  In other words, if $f(x, y)$ is a primitive form which is self-reciprocal or self-inversive, is it true that $f(x, y)$ is reduced?  This question is addressed in \cite{b-sh}.

\bigskip

\noindent \textbf{Acknowledgments:}  We would like to thank Lubjana Beshaj for  helpful conversations and explaining to us the reduction theory of self-inversive and self-reciprocal forms.  

\section{Self-inversive polynomials}\label{sect-1}

\def\Sl{\mbox{SL}}
\def\im{\mbox{Img }}
\def\a{\alpha}
\def\b{\beta}
\def\g{\gamma}
\def\d{\delta}
\def\psl{\mbox{PSL}}


Let $\P^1$ be the Riemann sphere and $ \gl_2(\C)$  the group of $2 \times 2$ matrices with entries in $\C$. Then $\gl_2(\C)$ acts on $\P^1$  by linear fractional transformations.     This action  is a transitive action, i.e. has only one orbit. 
Consider now the action of $\Sl_2(\R)$ on the Riemann sphere. This action   is not transitive, because   for $M= \begin{pmatrix} \alpha &  \beta \\ \gamma &  \delta  \end{pmatrix} \in \Sl_2(\R)$ we have
\[ \im \left(M z\right) = \frac{(\alpha \delta -\beta \gamma )}{|\gamma z+\delta|^2}  \im  z. \]
Hence, $z$ and $Mz$ have the same sign of imaginary part when $\det(M)=1$.  
The action of  $\Sl_2(\R)$  on $\P^1$ has three orbits,  namely $\R \cup \infty$, the upper half plane, and the lower-half plane.
Let  $\H_2$ be the complex upper half plane, i.e. 
\[\H_2 =\left \{ z = x+ i y \in \C \, \left | \frac{}{}\right. \, \im (z) >0\right \} \subset \C. \]
The group $\Sl_2(\R)$ preserves $\H_2$ and acts transitively on it, since for $g \in \Sl_2(\R)$ and $z \in \H_2$ we have
\[\im  (gz) = \frac{\im  z}{|\g z +\d|^2} > 0 \]
The modular group $\Gamma =\Sl_2(\Z) /\{\pm I\}$ also acts   on $\H_2$.  This action has a  fundamental domain  $\mathcal F$ 
\[\mathcal F= \left\{\frac{}{} z \in \mathcal H_2  \left| \frac{}{} \right.    \, |z|^2 \geq 1  \, \text{ and } \, |Re(z)| \leq 1/2 \right\}\]
%

Consider now all binary quadratic forms with real coefficients.  A quadratic form $f \in \R[x, y]$ has two complex roots (conjugate of each other) if $f$ is positive definite.  Hence, we have a one to one correspondence between positive definite quadratic forms and points of $\H_2$.  For a given $f\in \R[x, y]$, let $\xi (f)$ denote the zero of $f$ in $\H_2$.  This is called \textit{zero map}.  The positive definite binary form $f$ has minimal coefficients if and only if $\xi (f) \in \mathcal F$; see \cite{reduction} for details. 

The group $SL_2 (\R)$ acts on the set of positive definite quadratic forms by linear changes of coordinates.  Moreover, the zero map $f \mapsto \xi (f)$ is equivariant under this action. In other words, $\xi (f^M) = \xi (f)^M$, for any $M\in \Sl_2(\R)$.  Hence, to \textit{reduce} a binary quadratic $f$ with integer coefficients we simply compute $\xi (f)$ and then determine $M \in \Gamma$ such that $\xi (f)^M \in \mathcal F$.  Then, the quadratic $f^M$ has minimal coefficients.  

This approach can be generalized  to higher degree forms $f \in \R[x, y]$. 
 Then $f(x, y)$ is a product of linear and quadratic factors over $\R$.  
In studying roots of $f(x, y)$ we are simply concerned with roots in the upper half plane $\H_2$.  The zero map can also be defined in this case, but its definition is much more technical. The interested reader can check \cite{reduction} or \cite{beshaj} for details. 

Hence the problem of finding a form equivalent to $f$ with minimal coefficients becomes equivalent to determine a matrix $M \in \Gamma$ such that $\xi (f)^M \in \mathcal F$.  
The generators of the modular group $\Gamma$ are the matrices 
\[  S = \begin{bmatrix}  0 & -1 \\ 1 & 0 \end{bmatrix}   \quad \text{and} \quad 
T=\begin{bmatrix}  1 & 1 \\ 0 & 1 \end{bmatrix}
\]
which correspond to transformations $z \to - \frac 1 z$ and $z \to z+1$. 
Next, we will see the geometry  of some   of these transformations which play an  important role in this process. 

Let $ \sigma(z) = \frac 1 z$ be the \textbf{reciprocal map } of the complex plane.  Then, 
\[  \sigma ( a+ b i ) = \frac 1 {|z|^2} \, (a - bi) \]
Hence, on the unit circle $U = \{ z \in \C,  \, |z| =1     \}$ the reciprocal map  becomes simply the complex conjugation.  From this we see that to the geometric inversion of the unit circle corresponds  the \textbf{inversion map} 
\[ \tau :   z \to \frac  1  {\z} \]
which sends points $z \in \H_2$ inside the unit circle $U$  to points in $z^\prime \in  \H_2$ with the same argument as $z$ and $|z| \cdot |z^\prime| =1$.  
It fixes   points on the unit circle $U$.  It is exactly this transformation together with $z\mapsto z+1$  which    we use to "move" points within $\H_2$ and bring them in the fundamental domain. 
We are interested in forms $f(x, y)$ which are fixed by this transformation.  Hence, we are interested in polynomials $f(z, 1)$  whose set of roots is fixed by $\tau (z)$. 


For a degree $n$ polynomial $f(z)\in \C[z]$, the   \textbf{inversive} of $f$  is called the function $f^\star (z) = z^n f \left( \frac 1 {\z} \right) $. A polynomial $f$  will be called \textbf{self-inversive} if $f=f^\star$. We can make this definition more precise. 

Let $p (z) \in \C [z]$ such that 
\begin{equation}\label{eq-1}
p(z) = \sum_{i=0}^n a_i z^i.
\end{equation}
Then, $p (z)$  is called \textbf{self-inversive} if its set of zeroes is fixed by the inversion map $\tau (z) = 1/ {\z}$. Thus, the set of roots is 
\[ \left\{ \alpha_1 , \dots , \alpha_n,  \frac 1 {\bar \alpha_1}, \dots , \frac 1 {\bar \alpha_n} \right\} \]
and then $p(z)$ is given by 
\begin{equation}
 p(z) = a_n \, \prod_{i=1}^s  \, \left(  z^2 - \left( \alpha_i + \frac 1 {\bar \alpha_i} \right) z + \frac {\alpha_i} {\bar \alpha_i}  \right), 
\end{equation} 
Let us denote by $\p (z)$   the  \textit{conjugate polynomial}  of $p(z)$, namely   
\[ \p (z) := \sum_{i=0}^n {\bar a}_i z^i .\]
 Then, we have the following; see \cite{ohara}. 
\begin{lem}\label{lem1}
If $p(z)$ be given as in Eq.~\eqref{eq-1}.  The following are equivalent:

\begin{enumerate}
\item  $p(z)$ is self-inversive

\item  For every $z\in \C\setminus \{ 0\}$,  
\[   \aa_n \, p(z) = a_0 \, z^n \, \p \left( \frac 1 z  \right) \]
\item For every $z\in \C\setminus \{ 0\}$
\[ p (z) = w \cdot  z^n \cdot \p \left( \frac 1 z \right),\]
 where $|w|=1$.

\item For $j=0, 1, \dots , n$,
\[a_0 \aa_j = \aa_n a_{n-j}\]
\end{enumerate}

Moreover,  if $p(z)$ is self inversive 
then  
\begin{enumerate}
\item $ | a_i | = | a_{n-i}|$     for all $i=0, \dots , n$.    
\item    $\aa_n \left[ n \, p(z)  - z \, p^\prime (z) \right] = a_0 \, z^{n-1} \, \p^\prime \left( \frac 1 z \right) $, for each $z\in \C$
\item $\left| n \cdot \frac {p(z)} {z \cdot p^\prime (z)} - 1 \right| =1$, for each $z \in U$.   
\end{enumerate}
\end{lem}

Studying roots of the self-inversive polynomials is an old problem which has been studied by many authors.    A classical result due to Cohn states that a self-inversive polynomial has all its zeros on the unit circle if and only if all the zeros of its derivative lie in the closed unit disk.   

For $p(z) \in \C [z]$ we  let $||p||$ denote the \textit{maximum modulus} of $p(z)$ on the unit circle.  In \cite{ohara} it is proved the following

\begin{thm}
If $p (z) = \sum_{i=0}^n a_i z^j$, $a_n\neq 0$,  is a self-inversive polynomial which has all the zeroes on $|z|=1$, then
\[    |a_i| \leq \frac {||p||} 2 \]
for each $i \neq \frac  n 2$ and  $ |a_{n/2} | \leq \frac {\sqrt{2} } 2  ||p||$. 
\end{thm}

From the above theorem we can see that the middle coefficient is
special.  The \textit{middle coefficient conjecture} says that for $p(z)$ as in the above theorem, it is conjectured that 
\begin{equation}
\label{eqn:MCC}
| a_{n/2} | \leq ||p|| 
\end{equation}
If $n$ is even then the middle coefficient conjecture is true when $ | a_{n/2} | \leq 2 |a_n|$; see \cite[pg. 334]{ohara} for details.

The following theorem holds; see \cite{vieira}, \cite{lakatos} for details. 

\begin{thm}
Let $p \in \C [x]$ be a degree $n$ self-inversive polynomial.  If 
\[ 
| a_{n-\lambda} | > \frac 1 2 \, {n \choose {n - 2\lambda} }  \,  \, \sum_{k=0, k\neq \lambda, k\neq n-\lambda}^n |a_k| 
\]
for some $\lambda < \frac n 2$,  then $p(z)$ has exactly $n-2\lambda$ non-real roots on the unit circle. 

If $n$ is even and $\lambda = \frac n 2$, then $p(z)$ has no roots on the unit circle if 
\[ 
| a_{n/2} | > 2 \sum_{k=0, k \neq n/2}^n  | a_k   | 
\]
\end{thm}
For a proof see \cite{vieira}. If $\lambda =0$ this correspond to a result of Lakatos and Losonczi \cite{lakatos} which says that a self-inversive polynomial with non-zero discriminant has all roots on the unit circle if 
\[ 
| a_n | \geq \frac 1 2 \, \sum_{k=1}^n |a_k| \, .
\]
%

There is a huge amount of literature on bounding the roots or the coefficients of polynomials or finding polynomials which have bounded coefficients.  Most of that work relates to Mahler measure and related works. There was another approach by Julia \cite{julia} which did not gain the attention it deserved.   Lately there are works of Cremona and Stoll in \cite{SC},  Beshaj \cite{beshaj, reduction}, and others who have extended   Julia's method and provide an algorithm of finding the polynomial (up to a coordinate change) with the smallest coefficients.  The first paragraph of this section eludes to that approach.

\subsection{Reciprocal polynomials} 
For a degree $n$ polynomial $f(z) \in \C[z]$, its \textbf{reciprocal} is called the polynomial $f^\times(z) = z^n \, f \left(  \frac 1 z \right) $.  A polynomial is called \textbf{self-reciprocal} or \textbf{palindromic} if $f=f^\times$ and it is called \textbf{anti-palindromic} if $f= - f^\times$.  

If $p(z) \in \C[z]$ be a polynomial such that its set of roots is fixed by reciprocal map $\sigma (z)$, say 
\[ S = \left\{ \alpha_1, \dots , \alpha_s,  \frac 1 {\alpha_1},    \dots , \frac 1 {\alpha_s}\right\}, \]
then $f(z)$ is   palindromic or antipalindromic polynomial.  Due to the properties of the binomial coefficients the polynomials $P(x) = (x + 1 )^n$ are palindromic for all positive integers $n$, while the polynomials $Q(x) = (x - 1 )^n$ are palindromic when $n$ is even and anti-palindromic when $n$ is odd. Also, cyclotomic polynomials are palindromic.

What if  we would like some kind of invariant of the reciprocal map $z \mapsto 1/z$?    Consider  the transformation 
\[ \alpha (z)= z + \frac 1 z\]
Obviously, $\alpha ( 1/z) = z$.  
When considered as a function $\alpha : \C \to \C$ this is a 2 to 1 map since both $z$ and $1/z$ go to the same point. Considered on each one of the three orbits of $SL_2 (\R)$ in $\C$ we have the following:  $\alpha$ sends  the upper half-plane $\H_2$ onto the complex plane  $\C$ except for $(-\infty, 2]$ and $[2, \infty)$ 
 which are doubly covered by $\R \setminus \{0\}$.  
We organize such actions in the following Lemma: 

\begin{lem}\label{lem-3}
For any  polynomial $p (z)  = \sum_{i=0}^n a_i z^i$    of   degree $n=2s$ the following are equivalent:

\begin{enumerate}
\item  The coefficients of $p(z)$   satisfy 
\[ a_i = a_{n-i}, \quad \textit{ for all } \, i =0, \dots n.\]
\item  There exists a polynomial $q(z)$ such that 
\[ p(z) = z^s \, q \left( \frac 1 z \right)\] 
\item    There exists  some polynomial $g(z)$ of degree $m \geq 1$ such that 
 \[ p(z) = z^m \cdot g \left(  z + \frac 1 z \right) \]
\end{enumerate}
\end{lem}

For a proof see \cite{joyner} among other papers.    Hence, any polynomial $f(z)$ satisfying any of the properties of the Lemma is self-reciprocal.  

Next we list some properties of palindromic and antipalindromic polynomials.  Their proofs are elementary and we skip the details. 

\begin{rem} Here are some general properties of palindromic and anti-palindromic polynomials:
\begin{enumerate}

\item  For any antipalindromic  polynomial $p (z)  = \sum_{i=0}^n a_i z^i$
\[ a_i = - a_{n-i}, \quad \textit{ for all } \, i =0, \dots n.\]

\item   For any polynomial $f$, the polynomial $f + f^\times$  is palindromic and the polynomial $f - f^\times$  is antipalindromic.

\item   The product of two palindromic or antipalindromic polynomials is palindromic.

\item   The product of a palindromic polynomial and an antipalindromic polynomial is antipalindromic.

\item    A palindromic polynomial of odd degree is a multiple of $x+1$ (it has -1 as a root) and its quotient by $x+1$ is also palindromic.

\item   An antipalindromic polynomial is a multiple of $x-1$ (it has 1 as a root) and its quotient by $x-1$ is palindromic.

\item   An antipalindromic polynomial of even degree is a multiple of $x^2 - 1$  (it has -1 and 1 as a roots) and its quotient by $x^2-1$ is palindromic. 
\end{enumerate}
\end{rem}

The following lemma shows an important correspondence among the pairs of roots $\left(  \a , \frac 1 \a \right)$ of $f(z) $ and real roots of $g \left( z + \frac 1 z \right)$.  Polynomials $f(z)$ which have all roots on the units circle correspond to $g \left( z + \frac 1 z \right)$ which have all real roots.  When homogenized the corresponding forms are called \textit{totally real forms} (cf. Section 4).

\begin{lem}
Let $f(z)= \sum_{i=0}^n $ be a palindromic polynomial and $g(z) \in \C[z]$ such that $f(z)= z^m g(z + 1/z)$.  Denote by $S_f$ the set of pairs of roots of $f(z)$ on $U$,
\[ S_f = \left\{     \left( \alpha, \frac 1 \alpha \right),  \, \textit{such that } \, | \alpha| =1 \, \textit{ and } \, f(\alpha ) = 0 \right\} \]
and by $S_g$ the set of roots of $g(z)$ in $[ -2, 2 ]$.   There is a one-to-one correspondence between $S_f$ and $S_g$.  
\end{lem}

\proof
The proof is rather elementary. If $|z|=1$ then $z= \cos \theta + i \sin \theta$, for some $\theta$.  Then, $\alpha (z) = 2 \cos \theta $ is in the interval $[-2, 2]$.  
Conversely, if $t \in [-2, 2]$ then $t= 2 \cos \theta$ for some $\theta$.  Hence, $t= z + 1/z$, where $z = \cos \theta \pm i \sin \theta$. 
\endproof

Notice that the inversion $z \mapsto 1/z$ induces an involution on the group of symmetries of a palindromic polynomial.  Hence, the Galois group of such polynomials is non-trivial. We will see in the next section how such involution among the roots of $f(x)$ induces automorphisms for algebraic curves with affine equation $y^n = f(x)$.

A polynomial $f(z)= \sum_{i=0}^n a_i z^i$ is called \textbf{quasi-palindromic} if 
\[ | a_i | = | a_{n-i}| ,\]
for all $i=0, \dots , n$.

The following Lemma will be used in the next section.  

\begin{lem}\label{product}
Let $f, g \in \C [x]$ with no common factor.  If $f$ and $g$ are self-inversive then $f g$ is a self-inversive. If $f$ and $g$ are quasi-palindromic, then $fg$ is quasi-palindromic.
\end{lem}

\proof
The proof is an immediate consequence of the definitions.  Since the set of roots of $f$ and $g$ contain all $z$ and $\frac 1 {\z}$  (resp. $z$ and $\pm \frac 1 {z}$ ), then so would contain their union, which is the set of roots of $fg$.
\qed

\begin{rem}
A polynomial with real coefficients all of whose complex roots lie on the unit circle in the complex plane (all the roots are unimodular) is either palindromic or antipalindromic
\end{rem}









\subsection{Self-reciprocal  polynomials over the reals}
\label{sec:SRPs}


%
Here is a basic fact about even degree self-reciprocal polynomials; see \cite{DH}, \S 2.1; see also \cite{L2}. 
%
%
The degree $d=2n$ polynomial $p(z)$  is self-reciprocal if and only if it can be written
\[
p(z)=z^n\cdot (a_n+a_{n+1}\cdot (z+z^{-1})+\dots +
a_{2n}\cdot (z^n+z^{-n})),
\]
if and only if it can be written
\begin{equation}
\label{eqn:prod}
p(z)=a_{2n}\cdot \prod_{k=1}^n (1-\alpha_kz+z^2),
\end{equation}
for some real $\alpha_k\in \R$.

Note that $g(z)=1-\alpha z+z^2$ has roots on the unit circle if and only if the roots are of the form $e^{\pm i\theta}$, for some $\theta$, in which case, 
$\alpha = 2\cos(\theta)$.

%
%
%
For the rest of this section we denote by $p(z)=\sum_{i=0}^n a_iz^i$ a degree $n$ self-reciprocal polynomial, where $n=2d$ or  $n= 2d+1$.
The answer to the following question is unknown at this time: 
for which increasing sequences $a_0<a_1<\dots a_{d}$ do the roots of the
corresponding self-reciprocal
polynomial, $p(z)=0$, lie on the unit circle $|z|=1$?

If $n=2d$, which $p(z)$ with   $a_0<a_1<\dots a_{d}$,  can be written as a product $\prod_{k=1}^d (1-2\cos(\theta_k)z+z^2)$?

It is clear that, in a product such as (\ref{eqn:prod}),
with all its roots on the unit circle so $-2\leq \alpha_k \leq 2$, we have
\begin{equation}
\label{eqn:symm-increase}
0<a_0\leq a_1 \leq \dots \leq a_n,\ \ a_{n-i}=a_{n+i},
\end{equation}
for all $i\in \{0,1,2,\dots,n\}$, provided
the collection $\alpha_j$s satisfy
\begin{equation}
\label{eqn:symm-increase2}
\alpha_k\leq -1.
\end{equation}

A self-reciprocal polynomial satisfying (\ref{eqn:symm-increase}) is called    \textit{symmetric increasing}.
Motivated by Problem \ref{conj:symm-increasing} below, we look for a bound which is more general than (\ref{eqn:symm-increase2}) and which also implies
the polynomial is symmetric increasing.  For instance, we observe that the following result can be used inductively to establish a generalization of
 (\ref{eqn:symm-increase2}).

\begin{lem}
Let $p(z)$ be as above. To multiply 
$p(z)$ by $1-\alpha x+x^2$ ($-2\leq \alpha\leq 2$), and still have the new coefficients satisfy 
a symmetric increasing condition such as in (\ref{eqn:symm-increase}), we require

\begin{equation}
\label{eqn:symm-incr-cond}
(a_i,a_{i+1},a_{i+2},a_{i+3})\cdot (1,-1-\alpha,1+\alpha,-1)\leq 0,
\end{equation}
for all $i\leq d$.
In particular, if $a_i=a$, $a_{i+1}=a+\epsilon_1$, $a_{i+2}=a+\epsilon_2$, 
$a_{i+3}=a+\epsilon_3$ then  (\ref{eqn:symm-incr-cond}) holds if
\[
\epsilon_2 \leq \frac{\epsilon_1+\epsilon_3}{2}.
\]
\end{lem}

\proof    This is verified simply by multiplying out $p(z)(1-\alpha x+x^2)$, so omitted.

\qed

The examples below illustrate how sensitive (\ref{eqn:symm-increase}) is to the size of the $\alpha_j$s.

\begin{exa}
We have
{\small{
\[
(1+1.05x+x^2)(1-0.28x+x^2)(1+1.25x+x^2)=
x^6 + 2.02x^5 + 3.6685x^4 +
3.67250x^3 + 3.6685x^2 + 2.02x + 1,
\]
}}
which satisfies (\ref{eqn:symm-increase}), but change the $0.28$ to $0.3$ and 
{\small{
\[
(1+1.05x+x^2)(1-0.30x+x^2)(1+1.25x+x^2)=
x^6 + 2x^5 + 3.6225x^4 + 3.60625x^3 + 3.6225x^2 + 2x + 1,
\]
}}
does not.
Similarly, we have
{\small{
\[
(1+1.05x+x^2)(1-0.3x+x^2)(1+1.25x+x^2)(1-0.6x+x^2)=
\]
\[
x^8 + 1.4x^7 + 3.4225x^6 +
3.43275x^5 + 5.08125x^4 + 3.43275x^3 +
3.4225x^2 + 1.4x + 1,
\]
}}
which satisfies (\ref{eqn:symm-increase}), but change the $0.6$ to $0.7$ and 

{\small{
\[
(1+1.05x+x^2)(1-0.3x+x^2)(1+1.25x+x^2)(1-0.7x+x^2)=
\]
\[
x^8 + 1.3x^7 + 3.2225x^6 +
3.0705x^5 + 4.720625x^4 + 3.0705x^3 +
3.2225x^2 + 1.3x + 1,
\]
}}
does not.

The polynomial
\[
(1+1.5x+x^2)(1+0.2x+x^2)(1+0.1x+x^2)
=
x^6 + 1.8x^5 + 3.47x^4 + 3.63x^3 + 3.47x^2 + 1.8x + 1
\]
satisfies (\ref{eqn:symm-increase}), as does
\[
(1+1.5x+x^2)(1+0.2x+x^2)(1+0.1x+x^2)(1-0.5x+x^2)=
\]
\[
x^8 + 1.3x^7 + 3.57x^6 + 3.695x^5 + 5.125x^4 + 3.695x^3 + 3.57x^2 + 1.3x + 1
\]
but  change the $0.5$ to $0.6$ and the product
\[
(1+1.5x+x^2)(1+0.2x+x^2)(1+0.1x+x^2)(1-0.6x+x^2)=
\]
\[
x^8 + 1.2x^7 + 3.39x^6 + 3.348x^5 + 4.762x^4 + 3.348x^3 + 3.39x^2 + 1.2x + 1
\]
does not. 

The polynomial

\[
(1+0.1x+x^2)(1+0.2x+x^2)(1+0.3x+x^2)(1+0.92x+x^2)=
\]
\[
x^8 + 1.52x^7 + 4.662x^6 +
4.6672x^5 + 7.32952x^4 + 4.6672x^3 +
4.662x^2 + 1.52x + 1
\]
satisfies (\ref{eqn:symm-increase}), but  change the $0.92$ to $0.91$ and 

\[
(1+0.1x+x^2)(1+0.2x+x^2)(1+0.3x+x^2)(1+0.91x+x^2)=
\]
\[
x^8 + 1.51x^7 + 4.656x^6 + 4.6361x^5 + 7.31746x^4 + 4.6361x^3 + 4.656x^2 + 1.51x + 1
\]
does not.

The above lemma holds, namely the condition
(\ref{eqn:symm-incr-cond}), because

\[
(1+0.1x+x^2)(1+0.2x+x^2)(1+0.3x+x^2)(1+0.92x+x^2)(1+0.999x+x^2)=
\]
\[
x^{10} + 2.519x^9 + 7.18048x^8 +
10.844538x^7 + 16.6540528x^6 +
\]
\[
+ 16.65659048x^5 +
16.6540528x^4 + 10.844538x^3 + 7.18048x^2 +
2.519x + 1
\]
satisfies (\ref{eqn:symm-increase}), but change the $0.999$ to $0.99$ and 

\[
(1+0.1x+x^2)(1+0.2x+x^2)(1+0.3x+x^2)(1+0.92x+x^2)(1+0.99x+x^2)=
\]
\[
x^{10} + 2.51x^9 + 7.1668x^8 +
10.80258x^7 + 16.612048x^6 +
\]
\[
+ 16.5906248x^5 +
16.612048x^4 + 10.80258x^3 + 7.1668x^2 +
2.51x + 1
\]
does not.
\end{exa}

\def\Aut{\mbox{Aut } }
\def\bAut{\overline{\mbox{Aut}}\, \, }
\def\X{\mathcal X}
\def\d{\delta}
\def\normal{\triangleleft}
\def\<{\langle}
\def\>{\rangle}
\newcommand\G{\bar G}

\section{Superelliptic curves and self-inversive polynomials}\label{sect-2}

The following theorem connects self-reciprocal polynomials with a very special class of algebraic curves, namely superelliptic curves.  We follow the definitions and notation as in \cite{super-1}.   

Fix an integer $g\ge2$.  Let $\X_g$ denote a genus $g$ generic planar  curve defined over an algebraically closed field $k$ of characteristic $p \geq 0$. We denote by $G$ the full automorphism group of $\X_g$.  Hence,  $G$ is a finite group. Denote by $K$ the function field of $\X_g$ and assume that the affine equation of $\X_g$ is given some polynomial in terms of $x$ and $y$. 

Let $H=\< \tau \>$ be a cyclic subgroup of $G$ such that $| H | = n$ and $H $ is in the center of $G$, where $n \geq 2$. Moreover, we assume that the quotient curve $\X_g / H$ has genus zero. 
The \textbf{reduced automorphism group of $\X_g$ with respect to $H$} is called the group  $\G \, := \, G/H$, see \cite{super-1}.  

Assume  $k(x)$ is the  genus zero subfield of $K$ fixed by $H$.   Hence, $[ K : k(x)]=n$. Then, the group  $\G$ is a subgroup of the group of automorphisms of a genus zero field.   Hence, $\G <  PGL_2(k)$ and $\G$ is finite.   It is a classical result that every finite subgroup of $PGL_2 (k)$  is  isomorphic to one of the following: $C_m $, $ D_m $, $A_4$, $S_4$, $A_5$, \emph{semidirect product  of elementary   Abelian  group  with  cyclic  group}, $ PSL(2,q)$ and $PGL(2,q)$.

The group $\G$ acts on $k(x)$ via the natural way. The fixed field of this action is a genus 0 field, say $k(z)$. Thus, $z$
is a degree $|\G| := m$ rational function in $x$, say $z=\phi(x)$.


\begin{lem}
Let $\X_g$ be a superelliptic curve of level $n$ with $| \Aut (\X_g)| > n$.  Then, $\X_g$ can be written as 
\[ y^n = f(x^s),   \quad \textit{or } \quad y^n = x f(x^s)  \]
for some $s > 1$.  
\end{lem}

The proof goes similar as for the hyperelliptic curves as in \cite{issac}.  Since below we display all equations of such curves in such form then the Lemma is obviously true.

Next we focus on studying the nature of the polynomial $f(x)$   and its connections to self-inversive polynomials. We are assuming that the curves are of characteristic zero, so the reduced automorphism group is cyclic, dihedral, $A_4$, $S_4$, or $A_5$. The list of equations, including the full group of automorphisms, the dimension of the loci, and the ramification of the corresponding covers can be taken from \cite{Sa2}.

\begin{thm}
If the  reduced automorphism group of a superelliptic curve $\X$  is nontrivial or not isomorphic to a cyclic group, then $\X$ can be written with the affine equation 
\[ y^n = f(x) \quad \textit{or} \quad y^n = x \cdot f(x) \]
where $f(x)$ is a palindromic or antipalindromic polynomial. 
%
If the reduced automorphism group is isomorphic to $A_5$,  then $f(x)$ is a quasi-palindromic plynomial. 
\end{thm}

\proof
If $\bAut (\X) $ is isomorphic to a dihedral group $D_{2m}$, then the equation of $\X_g$ can be written as in one of the following cases
\[ 
\begin{split}
y^n = &\ F(x):= \prod_{i=1}^\d(x^{2m}+\lambda_ix^m+1) \\
y^n = &\          (x^m-1)\cdot F(x), \\
y^n = &\   x\cdot F(x), \\
y^n = &\  (x^{2m}-1)\cdot F(x),\\
y^n = &\   x(x^m-1)\cdot F(x), \\
y^n = &\  x(x^{2m}-1)\cdot F(x),
\end{split}
\]
The polynomial $F(x)$ is palindromic from Lemma~\ref{lem-3}.  The polynomials   $x^m-1$ and $x^{2m}-1$ are antipalindromic.
From Lemma~\ref{product}   the products $(x^m -1 ) F(x)$ and $(x^{2m} -1 ) F(x)$ are antipalindromic.    Hence, if the reduced automorphism group of a superelliptic curve is isomorphic to a dihedral group then the equation of the curve can be written as $y^2 = f(x)$ or $y^2= x f(x)$, where $f(x)$ can be chosen to be a palindromic or antipalindromic polynomial. \\

If $\bAut (\X) $ is isomorphic to   $A_4$, then the equation of $\X_g$ can be written as in one of the following cases
\[
\begin{split}
y^n = &\ G(x) \\
y^n =   &\ (x^4+2i\sqrt{3}x^2+1)\cdot G(x),\\
y^n =   &\  (x^8+14x^4+1)\cdot G(x), \\
y^n =   &\ x(x^4-1)\cdot G(x),  \\
y^n =   &\ x(x^4-1)(x^4+2i\sqrt{3}x^2+1)\cdot G(x), \\
y^n =   &\ x(x^4-1)(x^8+14x^4+1)\cdot G(x), 
\end{split}
\]
where 
\[ G(x):= \prod_{i=1}^\delta(x^{12}-\lambda_ix^{10}-33x^8+2\lambda_ix^6-33x^4-\lambda_ix^2+1)  \]
Notice that every factor of $G (x)$ is palindromic, hence $G(x)$ is also palindromic from Lemma~\ref{product}.  The polynomials $x^4+ 2 i \sqrt{3} x^2 +1$ and 
$x^8 + 14 x^4 +1$ are palindromic and therefore $(x^4+ 2 i \sqrt{3} x^2 +1) \, G(x)$ and $(x^8 + 14 x^4 +1) \, G(x)$ are palindromic.  When multiplied by $x^4-1$ such polynomials become antipalindromic since $x^4-1$ is antipalindromic. 
So the equation of the curve can be written as $y^2 = f(x)$ or $y^2= x f(x)$, where $f(x)$ can be chosen to be a palindromic or antipalindromic polynomial. 

If $\bAut (\X) $ is isomorphic to   $S_4$, then the equation of $\X_g$ can be written as in one of the following cases
\[
\begin{split}
y^n = &\     M(x)\\
y^n = &\     \left( x^8+14x^4+1 \right)  \cdot M(x)\\
y^n = &\     x(x^4-1) \cdot M(x)\\
y^n = &\      \left( x^8+14x^4+1 \right)  \cdot x(x^4-1) \cdot M(x)\\
y^n = &\       \left( x^{12}-33x^8-33x^4+1  \right)\cdot M(x)\\
y^n = &\      \left( x^{12}-33x^8-33x^4+1  \right)  \cdot \left( x^8+14x^4+1 \right)  \cdot M(x)\\
y^n = &\      \left( x^{12}-33x^8-33x^4+1  \right) \cdot x(x^4-1) \cdot M(x)\\
y^n = &\      \left( x^{12}-33x^8-33x^4+1  \right) \cdot \left( x^8+14x^4+1 \right)  \cdot x(x^4-1)  M(x)\\
\end{split}
\]
where 
\begin{equation*}
\begin{split}
M (x)=  \, \prod_{i=1}^\d    &   \left(  x^{24}+\lambda_ix^{20}+(759-4\lambda_i)x^{16}+2(3\lambda_i+1228)x^{12}  +  (759-4\lambda_i)x^8      \right. \\
                              &  \left.    +\lambda_ix^4+1 \right) \\ 
\end{split}
\end{equation*}
Since every factor of $M(x)$ is palindromic, then $M(x)$ is palindromic.  By Lemma~\ref{product} we have that  the equation of the curve can be written as $y^2 = f(x)$ or $y^2= x f(x)$, where $f(x)$ can be chosen to be a palindromic or antipalindromic polynomial. The antipalindromic cases correspond exactly to the cases when $x^4-1$ appears as a factor. \\

Let $\bAut (\X) $ is isomorphic to  $A_5$.  This case is slightly different from the other cases due to the fact that now the reduced group has an element of order 5 and $f(x)$ will be written as a decomposition of $x^5$.  So the change of coordinates $x \mapsto -x$ will preserve the sign for odd powers and change it for even powers of $x$. 

Let $\Lambda (x)$, $Q(x)$, $\psi (x)$ be as follows
\begin{equation*}
\begin{split}
\Lambda (x)=  \,  \prod_{i=1}^\d & \, (  x^{60}+a_1 x^{55} + a_2 x^{50}+ a_3 x^{45} +a_4  x^{40}+a_5 x^{35} +a_6 x^{30} - a_5 x^{25} + a_4 x^{20} \\
& - a_3 x^{15}  + a_2  x^{10}  - a_1 x^5 +1  ) \\ 
a_1 & = \lambda_i - 684 \\
a_2 & = 55\lambda_i + 157434 \\
a_3 & = 1205\lambda_i-12527460 \\
a_4 & = 13090\lambda_i+77460495 \\
a_5 & = 69585\lambda_i-130689144 \\
a_6 & = 134761\lambda_i - 33211924 \\  
Q (x)  = & \,  x^{30}+522x^{25}-10005x^{20}-10005x^{10}-522x^5+1, \\
\psi (x)  =  &  \, x^4+2i \sqrt{3} x^2 +1 \\ 
\end{split}
\end{equation*}
Then, the equation of $\X_g$ can be written as in one of the following cases
\[
\begin{split}
y^n =   &\ \Lambda(x)\\
y^n =   &\ x(x^{10}+11x^5-1) \cdot \Lambda(x)\\
y^n =    &\ x (x^{20}-228x^{15}+494x^{10}+228x^5+1)(x^{10}+11x^5-1)\cdot \Lambda(x)\\
y^n =    &\ (x^{20}-228x^{15}+494x^{10}+228x^5+1)\cdot \Lambda(x)\\
y^n =    &\ Q (x) \cdot \Lambda (x)\\
y^n =     &\ x(x^{10}+11x^5-1).\psi(x)\cdot \Lambda(x)\\
y^n =    &\ (x^{20}-228x^{15}+494x^{10}+228x^5+1)\cdot \psi(x)\cdot\Lambda(x)\\
y^n =    &\ (x^{20}-228x^{15}+494x^{10}+228x^5+1)(x(x^{10}+11x^5-1))\cdot\psi(x)\cdot\Lambda(x)\\
\end{split}
\]
Notice that $\Lambda (x)$ is a quasi-palindromic  polynomial since all its factors are so.  So are $Q(x), \psi (x)$ and the other factors. 
 By Lemma~\ref{product} we can say that in this case he equation of the curve can be written as $y^2 = f(x)$ or $y^2= x f(x)$, where $f(x)$ can be chosen to be a quasi-palindromic polynomial.

This completes the proof of the theorem.
\qed

In \cite{hidalgo} it is shown that if the group $H$ is unique in $G$ and the reduced group $G/H$ is not cyclic or nontrivial, then the field of moduli is a field of definition for superelliptic curves. In \cite{beshaj-2016} and \cite{b-sh} it is explored the fact that most palindromic or self-inversive polynomials have minimal  coefficients.  So it is a natural question to investigate what is the relation between the minimal of definition of such curves, the minimal height as in \cite{beshaj}, and the palindromic polynomial $f(x)$.  

\def\J{\mathcal J}
\def\z{\varepsilon}

\section{Self-reciprocal polynomials and reduction theory}\label{sect-3}

Every stable binary form $f(z, y)$ of degree $n\geq 2$ correspond uniquely to a positive definite quadratic $\J_f$ called Julia quadratic; see \cite{beshaj}.  Since positive definite quadratics have a unique zero in the upper half plane $\H_2$, then we associate the zero of $\J_f$ to the binary form $f$.  This defines a map $\z$ from the set of degree $n$ binary forms to $\H_2$, which is called the zero map.  A binary form $f(z, y)$ is called \textit{reduced} if $\z (f) \in \mathcal F_2$. The size of the coefficients of a reduced binary form is bounded by its Julia invariant $\theta (f)$.  If $f$ is a reduced form, we say that $f$ has \textit{minimal coefficients}; see \cite{beshaj} for details.

There are no efficient ways to compute the Julia quadratic or the Julia invariant of a binary form of high degree (i.e. degree $> 6$).  Moreover, there is no known method to express the Julia invariant $\theta (f)$ in terms of the generators of the ring of invariants of the degree $n$ binary forms (i.e. transvections of the form). However, as discussed in \cite{beshaj} the case when $f$ is totally real is much easier.  A form is called \textit{totally real} if it splits over $\R$.


Let $f\in \C [z]$ be a degree $n\geq 2$ polynomial.  We denote by $f_*$ the corresponding form (homogenization of $f$) in $\C [z, y]$.      $GL_2 (\C )$ acts on the space of degree $n$ binary forms.  For a matrix $M \in GL_2 (\C)$ we denote by $f_*^M$ the action of $M$ on $f_*$.  By $f^M$ we denote $f_*^M (z, 1)$.  

\begin{lem}\label{lem-1} 
Let $f\in \C[z]$  and $M=\begin{bmatrix}   1 & -i \\ 1 & i \end{bmatrix}$.   Then,  
$f_*$ is a totally real binary form   if and only if $f^M$ has all roots in the unit circle.
\end{lem}

\proof
The proof is rather elementary.  The M\"obious transformation $h(z)=Mz$ maps $\H_2$ onto the open unit disk.  Moreover, it maps bijectively $U \setminus \{ 1\}$ to $\R$.  

\qed

For reduction of totally real forms see \cite{beshaj} and \cite{b-sh}.

\begin{thm}
Let $f(z)$ be a self-inversive polynomial.   Then   the following are equivalent:

i) all roots of $f(z)$ are on the unit circle

ii) all roots of its derivative $f^\prime (z)$ are on the unit disk

iii) $f_*^M$ is totally real form
\end{thm}

\proof  The equivalence of i) and iii) is the above Lemma.  The equivalence of i) and ii) is a result of Cohn. 
\qed

It is interesting to see how the reduction is performed in such case.   
From \cite{beshaj} we have a polynomial $G_f$ associated to $f$.
The  Julia quadratic $J_f$ is the only quadratic factor of $G_f$ when factored over $\R$.  Moreover, Beshaj \cite{beshaj} has proved that $G_f$ is very similar to a self-inversive polynomial. We describe briefly below


Let $f$ be a generic totally real form given by 
\[ f(x, y) = a_n x^n + a_{n-1} x^{n-1} y + \cdots + a_1 x y^{n-1} + a_0 y^n \]
where $a_0, \dots , a_n$ are transcendentals. 
Identify $a_0, \dots, a_n$ respectively with $1, \dots, n+1$. Then  the symmetric group  $S_{n+1}$ acts on $\R[a_0, \dots a_n ] [x, y]$ by permuting $a_0, \dots , a_n$.  For any permutation $\tau \in S_{n+1}$ and $f \in \R[a_0, \dots a_n ] [x, y]$ we denote by $\tau (f) = f^\tau$. Then 
\[ f^\tau (x, y) = \tau (a_n) \, x^n + \tau (a_{n-1} ) \, x^{n-1} y + \dots + \tau (a_1) \, x y^{n-1} + \tau (a_0) \, y^n. \]
Define $G(x, y)$    as follows
\begin{equation}\label{G-covariant}
G(x, y)= \frac  {    x \cdot f_x (- f_y (x, y), f_x (x, y) )+ y \cdot f_y (- f_y (x, y), f_x (x, y) )    }    {n \, f(x, y) }.
\end{equation}
Notice that since $f$ is totally real, then $f\in \R[x, y]$.  Therefore, $G \in \R [x, y]$.
Note also that,  for  $\sigma \in S_{n+1}$ we have   an involution  
\[
\sigma=
\left\{
\begin{split}
& (1,n+1)(2, n) \cdots  \left( \frac n 2 ,  \frac n 2 +2 \right),                                    &   \text{ if   $n$ is even }  \\
&  (1, n+1)(2, n) \cdots  \left(    \frac {n+1} 2,  \frac {n+3} 2   \right),     &   \text{ if   $n$ is odd. }
\end{split}
\right.
\]
Next result describes the properties of $G (x, y)$. 
\begin{thm}[Beshaj]  \label{thm1-totally-real}    The polynomial $G(x, y)$ satisfies the following \\

i) $G(x, y)$ is a covariant of $f$ of degree $(n-1)$ and order $(n-1)(n-2)$.  

ii) $G(x, y)$ has a unique quadratic factor  over $\R$, which is the Julia quadratic   $\J_f$.  

iii)  $G^\sigma (x, y)= (-1)^{n-1} \, G (x, y)$.  Moreover, if $ G_f = \sum_{i=1}^d   g_i \, x^i y^{d-i}$,   then 
   \[g_i^\sigma = (-1)^{n-1} \, g_{d-i},\]
    for all $i=0, ..., d$.
\end{thm}

Then we have the interesting connection between real forms and self-inversive polynomials. 

\begin{thm}
If $f$ is a palindromic real form then $G_f (x, y)$ is self-inversive.  If $f$ is of odd degree then $G_f$ is palindromic.  
\end{thm}

\proof
If $f$ is palindromic, then from Lemma 3, i) we have that $a_i = a_{n-i}$ for all $i=0, \dots , n$.  That means that $\sigma$ fixes all coefficients of $f$.  
Hence, $g_i^\sigma= (-1)^{n-1} \, g_i$ for all $i=0, \dots , d$, where $d= \deg G_f$.  Thus, $G_f$ is self-inversive.   If $n$ is odd, then 
$g_i^\sigma=   g_i$.  Hence, $G_f$ is palindromic.

\qed

We know that $G_f$ has exactly two non-real roots, namely $\z (f)$ and its conjugate.  
Consider now $G_f^M$.  Then all real roots of $G_f$ will go to roots on the unit circle of $G_f^M$ and the two non-real roots $\z (f)$ and its conjugate $\overline {\z (f)}$ go inside the unit disk as roots of $G_f^M$.  

\begin{figure}[h] 
   \centering
   \includegraphics[width=3in]{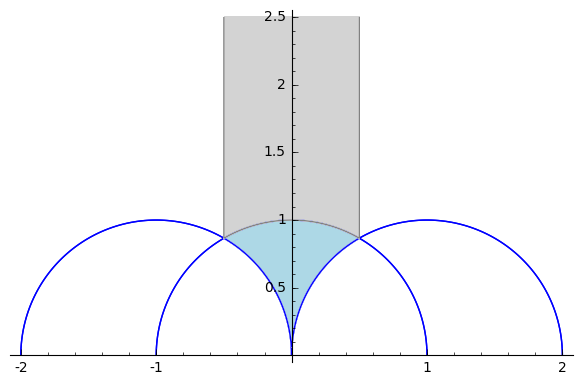} 
   \caption{The region $\mathcal T$}
   \label{fig1}
\end{figure}
\begin{lem} Let $f$ be  a self-inversive polynomial with all roots in the unit circle $U$, $f_*$ its homogenization,  $\mathcal T$  be the region in the complex plane given by
\[ \mathcal T = \{ z= a + b i \, |  \,  a^2 - 2a +b^2 \geq 0,   \, a^2 +2a+b^2 \geq 0  \} , \]
and $M=\begin{bmatrix}   1 & -i \\ 1 & i \end{bmatrix}$.  
If $ \z (f_*)^M \in \mathcal T$ or $ \z (f_*)^M \in \mathcal F_2$, then $f^M$ has minimal coefficients.  
\end{lem}

\proof
From Lem.~\ref{lem-1} we have that $f_*^M$ is a totally real form.  Then $\z (f_*^M)$ is the image of the zero map in the upper half plane $\H_2$.  

If $ \z (f_*)^M \in \mathcal F_2$ then $f_*^M$ is reduced and we are done.  If $ \z (f_*)^M \in \mathcal T$
then let  $S=\begin{bmatrix} 0 & 1 \\ 1 & 0 \end{bmatrix}$ and compute   $\z (f_*)^{MS}$. Let $\z (f_*)^M= a+bi$. 
Then 
\[ \z (f_*)^{MS} = \frac 1 {a+bi} = \frac a {a^2+b^2} - \frac b {a^2+b^2} i \]
Hence, $| \z (f_*)^{MS} | \geq 1$ and   
\[ - \frac 1 2 \leq \frac a {a^2 + b^2 } \leq \frac 1 2 \]
Hence, $ \z (f_*)^M \in \mathcal F_2$ 
However, the height of $f_*^M$ does not change under the transformation $S$.  Hence,   $f_*^M$ has minimal coefficients. Thus, in both cases $f^M$ has minimal coefficients. 
\qed

The region $\mathcal T$ is the blue colored region in Fig.~\ref{fig1} and the grey area is the fundamental domain.  


\section{Self-reciprocal polynomials and  codes}\label{sect-4}

The goal of this section is to show how self-reciprocal polynomials are connected to other areas of mathematics, namely whether extremal formal weight enumerators for codes satisfy the Riemann hypothesis.     We will follow the  setup of \cite{e-sh}.

For $d\leq n$, denote the weight enumerator of an MDS code $C$ over $\F = GF(q)$
of length $n$ and minimum distance $d$ by $M_{n,d}(x,y)$. The dual $C^{\perp}$ is also an MDS code of length $n$ and minimum distance $d^{\perp}=n+2-d$.  Therefore, for $d\geq 2$, the weight enumerator of $C^{\perp}$ is $M_{n,n+2-d}(x,y)$. Let $M_{n,n+1}=x^n$. The MDS code with weight enumerator $M_{n,1}$ has dimension $n-d+1=n-1+1=n$, hence $C=\F^n_q$. It is easy to see that $M_{n,n+1}$ is the MacWilliams transform,
$(x,y)\mapsto (\frac{x+(q-1)y}{\sqrt{q}},\frac{x-y}{\sqrt{q}})$, 
of $M_{n,1}$. We may think of $M_{n,1}$ as the weight enumerator of 
the zero code.

The set $\{M_{n,1}, M_{n,2},\ldots,M_{n,n-1}, M_{n,n+1}\}$ is a basis for the vector space of homogeneous polynomials of degree $n$ in $x,y$. Furthermore, this set is closed under the MacWilliams transform; see \cite{e-sh} for details. 

If $C$ is an $[n,k,d]_q$-code, then one can easily see that

\[
A_C(x,y)=\sum_{i=d}^{n+1}a_{i-d}M_{n,i}=a_0M_{n,d}+\ldots+a_{n+1-d}M_{n,n+1},
\]
for some integers $a_i$ as in \S 4.4.2 in \cite{JK}.
The zeta polynomial of $C$ is defined as 
\[P(T):=a_0+a_1T+ \dots +a_{n-d+1}T^{n+1-d}.\]
The zeta polynomial $P(T)$ of an $[n,k,d]_q$-code $C$ determines uniquely the weight enumerator of $C$. The degree of $P(T)$ is at most $n-d+1$.
 The quotient 
\[ Z(t)=\frac{P(T)}{(1-T)(1-qT)}\]
is called \textbf{the zeta function} of the linear code $C$. 
The zeta function of an MDS code 
\[ \frac{1}{(1-T)(1-qT)}=\sum_{j=0}^{\infty}\frac{q^{j+1}-1}{q-1}T^j\]
is the rational zeta function over $\F_q$; see \cite[Cor. 1]{e-sh}. 
%
%
%
%
%
%
Formally self-dual codes lead to self-reciprocal polynomials.  The proof of the following Proposition can be found in \cite{e-sh}. 

\begin{prop} \label{self reciprocal} 
If $P(T)$ is the zeta polynomial of a formally self-dual code, then $P(T/{\sqrt q})$ is a self-reciprocal polynomial.
\end{prop}

\subsection{Riemann zeta function versus zeta function  for self-dual codes}

From \cite{e-sh} we have  that for a self-dual code $C$,  
\[ Z(T)=q^{g-1}T^{2g-2}Z(1/qT),\] which for  $z(T):=T^{1-g}Z(T)$, 
may be written as 
\[ z(T)=z(1/qT).\]
Now let 
\[ \zeta_C(s):=Z(q^{-s}) \, \,   \text{and}   \, \,     \xi_C(s):=z(q^{-s}).\]
We obtain 
\[ \xi_C(s)=\xi_C(1-s),\]
which is the same symmetry equation is analogous to the functional equation for
the Riemann zeta function. We note that $\zeta(s)$ and $\xi(s)$ have the same zeros.

The zeroes of the zeta function of a linear code $C$ are useful in understanding possible values of its minimum distance $d$.

Let $C$ be a linear code with weight distribution vector $(A_0,A_1,\ldots,A_n).$ Let $\alpha_1,\ldots, \alpha_r$ be the zeros of the zeta polynomial $P(T)$ of $C$ Then
\[ d=q-\sum_{i}{\alpha_i^{-1}}-\frac{A_{d+1}}{A_d}\frac{d+1}{n-d}. \]
In particular,
\[ d\leq q- \sum_ {i}{\alpha_i^{-1}}; \]
see \cite{e-sh} for details.

A self-dual code $C$ is said to satisfy {\it Riemann hypothesis} 
if the real part of any zero of $\zeta_C(s)$ is $1/2$, or equivalently, 
the zeros of the zeta polynomial $P_C(T)$ lie on the circle 
$|T|=1/\sqrt q$, or equivalently, the roots of the self-reciprocal 
polynomial (see Proposition \ref{self reciprocal} above) 
$P_C(T/\sqrt q)$ lie on the unit circle.

While Riemann hypothesis is satisfied for curves over finite 
fields, in general it does not hold for linear codes. A result 
that generates many counterexamples may be found in \cite{JK}. 
There is a family of self-dual codes that satisfy the Riemann 
hypothesis which we are about to discuss. The theory involved 
in this description holds in more generality than linear codes 
and their weight enumerators. 

\subsection{Virtual weight enumerators}
A homogeneous polynomial
\[ F(x,y)=x^n+\sum_{i=1}^{n}{f_ix^{n-i}y^i}\]
with complex coefficients is called a \textit{virtual weight enumerator}. The set 
\[ \{0\}\cup \{i:f_i\neq 0\}\]
 is called its \textit{support}. If 
\begin{equation}\label{eq-F}
 F(x,y)=x^n+\sum_{i=d}^n f_ix^{n-i}y^i, 
\end{equation}
with $f_d\neq 0,$ then $n$ is called the \textit{length} and $d$ is called the \textit{minimum distance} of $F(x,y)$. 

Let $C$ be a self-dual linear $[n,k,d]$-code. Recall that $n$ is even, $k=n/2$ and its weight enumerator satisfies MacWilliams' Identity.  A virtual generalization of $A_C(x,y)$ is straightforward. A virtual weight enumerator $F(x,y)$ of even degree that is a solution to MacWilliams' Identity
\begin{equation}
\label{eqn:VSD}
F(x,y)=F \left(\frac{x+(q-1)y}{\sqrt{q}},\frac{x-y}{\sqrt{q}} \right),
\end{equation}
is called \textit{virtually self dual} over $\F_q$ with \textit{genus} $\gamma(F)=n/2+1-d$. Although a virtual weight enumerator in general does not depend on a prime power $q$, a virtually self-dual weight enumerator does. 

\begin{prob}
Find the conditions under which a (self-dual) virtual weight enumerator with positive integer coefficients arises from a (self-dual) linear code. 
\end{prob}

The zeta polynomial and the zeta function of a virtual weight enumerator are defined as in the case of codes.

\begin{prop} [\cite{Ch}]
Let $F(x,y)$ be a virtual weight enumerator of length $n$ and minimum distance $d$. Then, there exists a unique function $P_F(T)$ of degree at most $n-d$ which satisfies the following
\[\frac{(y(1-T)+xT)^n}{(1-T)(1-qT)}P_F(T)=\ldots +\frac{F(x,y)-x^n}{q-1}T^{n-d} +\ldots\]
\end{prop}

\noindent The polynomial $P_F(T)$ and the function 
\[ Z_F(T):=\frac{P(T)}{(1-T)(1-qT)},\]
 are called respectively \textit{the zeta polynomial and the zeta function of the virtual weight enumerator} $F(x,y)$.

A virtual self-dual weight enumerator \textit{satisfies the Riemann hypothesis} if the zeroes of its zeta polynomial $P_F(T)$ lie on the circle $|T|=1/\sqrt{q}$.
There is a family of virtual self-dual weight enumerators that satisfy Riemann hypothesis. It consists of enumerators that have certain divisibility properties.

Let $b>1$ be an integer. If supp$(F)\subset b \Z$, then $F$ is called $b$-\textit{divisible}.
%
Let $F$ given  by  Eq.~\eqref{eq-F}
%
%
 be a $b$-divisible, virtually self-dual weight enumerator over $\F_q$. Then   $F(x,y)$   is called
\begin{description}
\item [Type I] if $q=b=2,~2|n$.
\item [Type II] if $q=2, b=4, 8|n$.  
\item [Type III] if $q=b=3, 4|n$.
\item [Type IV] if $q=4, b=2, 2|n$.
\end{description}

\noindent Then we have the following theorem:
\begin{thm}[Mallows-Sloane-Duursma]  \label{bound}
If $F(x,y)$ is a $b$-divisible self-dual virtual enumerator with length $n$ and minimum distance $d$, then
\[
d \leq 
\left\{ 
\begin{split}
& 2 \left[  \frac{n}{8}\right] +2,  \quad & \text{ if  F is Type I},\\
& 4 \left[ \frac{n}{24}\right]  +4, \quad & \text{ if F is Type II}, \\
& 3  \left[ \frac{n}{12}\right] +3, \quad &  \text{ if F is Type III}, \\
& 2  \left[ \frac{n}{6}\right] +2, \quad &  \text{ if F is Type IV}.\\
\end{split}
\right.
\]
\end{thm}

A virtually self-dual weight enumerator $F(x,y)$ is called \textit{extremal} if the bound in Theorem \ref{bound} holds with equality. 
%
A linear code $C$ is called $b$-divisible, extremal, Type I, II, II, IV if and only if  its weight enumerator has the corresponding property.

\noindent The zeta functions of all extremal virtually self-dual weight enumerators are known; see \cite{D3}. 
The following result can be found in \cite{D3}. 

\begin{prop} 
All extremal type IV virtual weight enumerators satisfy the Riemann hypothesis.
\end{prop}

\noindent For all other extremal enumerators, Duursma has suggested the following conjecture in \cite{D4}.

\begin{prob} 
Prove that any extremal virtual self-dual weight enumerators of type I-III satisfies the Riemann hypothesis.
\end{prob}


Let $F$ denote a weight enumerator as in (\ref{eqn:VSD}) and $P_F(T)$ the associated zeta polynomial.
Let $p_F(T)=P_F(T/\sqrt{q})$ denote the normalized zeta polynomial. Numerous computations suggest the
following result.

\begin{prob} 
\label{conj:symm-increasing}
If $F$ is an extremal weight enumerator of Type I, II, II, IV then the normalized zeta polynomial
is symmetric increasing. In fact, using the notation of (\ref{eqn:symm-increase}), if
if $a_i=a$, $a_{i+1}=a+\epsilon_1$, $a_{i+2}=a+\epsilon_2$, 
$a_{i+3}=a+\epsilon_3$ then    $\epsilon_2 \leq \frac{\epsilon_1+\epsilon_3}{2}$.
\end{prob}




\bibliographystyle{amsplain}

\bibliography{mybib}{}


\end{document}